\font\titlefont=cmr12
\font\sectionfont=cmbx10 at 11 pt

\magnification=\magstep1
\baselineskip=15pt
\mathsurround=1pt
\abovedisplayskip=8pt plus 3pt minus5pt
\belowdisplayskip=8pt plus 3pt minus5pt
\fontdimen16\textfont2=2.5pt  
\fontdimen17\textfont2=2.5pt  
\fontdimen14\textfont2=4.5pt  
\fontdimen13\textfont2=4.5pt  

\def\pf #1. {\noindent {\it #1}.\enskip}

\def\qed{\quad\hfill \rlap{$\sqcup$}$\sqcap$\par\medskip} 

\def\:{\,\colon} 

\def\bdy{\partial}

\outer\def\section#1\par{\vskip0pt plus.05\vsize\penalty-250\vskip0pt plus-.05\vsize 
\vskip1cm
\leftline{\sectionfont#1}\medskip}


\vbox{}\vskip 6mm
\centerline{\titlefont The Cyclic Cycle Complex of a Surface}

\vskip3mm
\centerline{Allen Hatcher}

\vskip12mm

A recent paper [BBM] by Bestvina, Bux, and Margalit contains a construction of a cell complex that gives a combinatorial model for the collection of isotopy classes of oriented curve systems within a given homology class in a closed orientable surface.  A key property of the complex is that it is contractible.  Here we will construct a different complex, this one a simplicial complex rather than a cell complex, with similar properties.

Here are the definitions.  By a {\it cycle\/} in a closed oriented surface $S$ we mean a nonempty collection of finitely many disjoint oriented smooth simple closed curves. A cycle $c$ is {\it reduced\/} if no subcycle of $c$ is the oriented boundary of one of the complementary regions of $c$ in $S$ (using either orientation of the region).  In particular, a reduced cycle contains no curves that bound disks in $S$, and no pairs of circles that are parallel but oppositely oriented. 

Define the cycle complex $C(S)$ to be the simplicial complex having as its vertices the isotopy classes of reduced cycles in $S$, where a set of $k+1$ distinct vertices spans a $k$-simplex if these vertices are represented by disjoint cycles $c_0,\cdots,c_k$ that cut $S$ into $k+1$ cobordisms $C_0,\cdots,C_k$ such that the oriented boundary of $C_i$ is $c_{i+1} - c_i$, subscripts being taken modulo $k+1$, where the orientation of $C_i$ is induced from the given orientation of $S$ and $-c_i$ denotes $c_i$ with the opposite orientation.  The cobordisms $C_i$ need not be connected.  The faces of a $k$-simplex are obtained by deleting a cycle and combining the two adjacent cobordisms into a single cobordism.  One can think of a $k$-simplex of $C(S)$ as a ``cycle of cycles". The ordering of the cycles $c_0,\cdots,c_k$ in a $k$-simplex is determined up to cyclic permutation.

Cycles that span a simplex represent the same element of $H_1(S)$ since they are cobordant.  Thus we have a well-defined map $\pi_0C(S)\to H_1(S)$. This has image the nonzero elements of $H_1(S)$ since on the one hand, every cycle representing a nonzero homology class contains a reduced subcycle representing the same class (subcycles of the type excluded by the definition of reduced can be discarded one by one until a reduced subcycle remains), and on the other hand, it is an elementary fact, left as an exercise, that a cycle that represents zero in $H_1(S)$ is not reduced.  

For a nonzero class $x\in H_1(S)$ let $C_x(S)$ be the subcomplex of $C(S)$ spanned by vertices representing $x$, so $C_x(S)$ is a union of components of $C(S)$.  In fact it is a single component, according to our main result:

\proclaim Theorem. $C_x(S)$ is contractible for each nonzero $x\in H_1(S)$.

When $S$ is the torus there is a unique reduced cycle in each nonzero homology class (unique up to isotopy, that is).  For a primitive class this cycle is a single curve, and for a class that is $p$ times a primitive class, $p>1$, the cycle consists of $p$ parallel copies of the curve representing the primitive class, with parallel orientations.  Thus $C_x(S)$ is a single point for each nonzero $x$.

The mapping class group of $S$ acts transitively on primitive homology classes, so the complexes $C_x(S)$ for primitive $x$ are all isomorphic.  This holds also for classes $x$ that are $p$ times a primitive class, for fixed $p$. 

\proclaim Proposition. If $S$ has genus $g > 1$ then $C_x(S)$ has dimension $2g-3$.  

\pf Proof. First we use Euler characteristic to show $2g-3$ is an upper bound on the dimension.  This amounts to showing that the cycle of cobordisms $C_i$ associated to a cycle of cycles $c_i$ defining a simplex of $C(S)$ can consist of at most $2g-2$ cobordisms.  Since we are dealing with reduced cycles, no component of a cobordism $C_i$ can be a disk, but some components can be annuli, provided that these annuli go from $c_i$ to $c_{i+1}$.  However, not all components of $C_i$ can be such annuli, otherwise $c_i$ would be isotopic to $c_{i+1}$, contrary to the definition of $C(S)$. Thus $\chi(C_i)\le -1$ for each $i$.  As $\sum_i\chi(C_i) = \chi(S) = 2-2g$ this implies there are at most $2g-2$ cobordisms in the given cycle of cycles.

Every simplex of $C(S)$ is contained in a simplex of dimension $2g-3$.  To see this, note first that every component of a cobordism $C_i$ of a simplex of $C(S)$ must have at least one boundary circle in each of $c_i$ and $c_{i+1}$.  Such a component which is not an annulus can be subdivided as a composition of pair-of-pants cobordisms. By adding extra annulus components where needed, this allows each $C_i$ to be subdivided as a chain of cobordisms each having $\chi = -1$. The total number of such cobordisms is then $2g-2$. \qed

Before proving the theorem we make a few remarks on the criterion for a cycle to be reduced.  A cobordism $C$ has its boundary partitioned into two collections of circles $\bdy_+ C$ and $\bdy_- C$. If one of these is empty we call $C$ a {\it null cobordism}.  Thus the condition for a cycle $c$ to be reduced can be restated as saying that the complementary cobordism $C$ has no components that are null cobordisms.  An elementary observation is that if a cobordism $C$ is the composite of two cobordisms $C'$ and $C''$, with $\bdy_+ C'=\bdy_- C''$, then if $C'$ and $C''$ have no components that are null cobordisms, the same is true for $C$.  As a consequence, given a cycle of cycles $c_0,\cdots,c_k$ with corresponding cobordisms $C_0,\cdots,C_k$, then the $c_i$'s are reduced if and only if no $C_i$ has a component that is a null cobordism.

\medskip
\pf First proof of the Theorem, for $x$ primitive.  This will be by a simple surgery argument in the spirit of [H].  Let $b$ be an oriented simple closed curve in the homology class of $x$.  This exists since $x$ is primitive.  We will construct a deformation retraction of $C_x(S)$ onto the star of $b$.  Since the star is contractible, this will prove the result.

A point in a simplex $[c_0,\cdots,c_k]$ of $C_x(S)$ can be regarded as a linear combination $c(t)=\sum_i t_i c_i$ with $\sum_i t_i=1$ and $t=(t_0,\cdots,t_k)$ where the $t_i$'s are the barycentric coordinates of the point. We can assume the $c_i$'s have been isotoped to intersect $b$ in the minimum number of points.  Such a minimal position is unique up to isotopy that preserves minimal position; this will guarantee that our constructions are independent of the choice of representative of an isotopy class.  Let $c=\cup_i c_i$.  Assuming that $c\cap b$ is nonempty, we will do a surgery process to eliminate all the points of $c\cap b$.  The process will give a path $c_u(t)$ in $C_x(S)$ starting at the given point $c(t)=c_0(t)$ and ending at a point in the star of $b$.

The orientations of the $c_i$'s and $S$ induce normal orientations of the $c_i$'s, with the normal orientation of $c_i$ pointing into the cobordism $C_i$.  These normal orientations give the points of $c\cap b$ normal orientations in $b$.  Since each $c_i$ represents the same homology class as $b$, the algebraic intersection numbers $c_i\cdot b$ are zero, hence also $c\cdot b=0$.  This means that there are as many points of $c\cap b$ with a clockwise orientation in $b$ as there are points with a counterclockwise orientation.  Hence for at least one of the arcs $\alpha_j$ into which $c$ cuts $b$, both endpoints of $\alpha_j$ have normal vectors pointing into $\alpha_j$.  We call such an $\alpha_j$ {\it innermost}.  

If $\alpha_j$ is innermost, it lies in some $C_i$ with both its endpoints on $c_i$.  We can use all the innermost $\alpha_j$'s in $C_i$ to surger $c_i$ to a new cycle $c'_i$ in the interior of $C_i$ which divides $C_i$ into two cobordisms $C'_i$ and $C''_i$ where $C'_i$ is the union of a collar neighborhood of $c_i$ and a neighborhood of the innermost $\alpha_j$'s in $C_i$.  If $c'_i$ is not reduced, there will be components of $C''_i$ with their boundary contained in $c'_i$.  These components can contain no arcs $\alpha_j$, as these arcs would have to be among the innermost $\alpha_j$'s that are already being used for the surgery. Adding these components of $C''_i$ to $C'_i$ and deleting their boundaries from $c'_i$ produces a new $c'_i$ that is reduced, since the new $C'_i$ and $C''_i$ have no components that are null cobordisms. In terms of barycentric coordinates, we let the coefficient $t_i$ of $c_i$ decrease to $0$ at unit speed while increasing the coefficient $t'_i$ of $c'_i$ at the same rate.  If there are other innermost arcs $\alpha_j$ in other $C_i$'s we treat them in the same way at the same time.  Once a coefficient $t_i$ decreases to $0$, the corresponding $c_i$ is deleted and we continue with the remaining innermost arcs.  

If we did not have to deal with the coefficients $t_i$, it would be clear that finitely many iterations of the step of surgering along all innermost $\alpha_j$'s simultaneously would produce a cycle of cycles disjoint from $b$.  To see that this still holds in the presence of the $t_i$'s it is helpful to replace each $c_i$ by a small bicollar neighborhood of itself in $S$, thinking of $t_i$ as the thickness of this bicollar, so that as $t_i$ shrinks to $0$ the bicollar shrinks to just $c_i$ itself and disappears.  Surgery on a cycle $c_i$ can then be regarded as a continuous operation, cutting gradually into the bicollar on $c_i$ to produce a progressively thinner bicollar of $c_i$ and a progressively thicker bicollar of $c'_i$.  This is the analog for cycles of what is done for arcs in [H].  The bicollars of $c$ intersect $b$ in a collection of disjoint oriented intervals about the points of $c\cap b$.  From these intervals we can produce a finite tree $T=T(c)$ by first collapsing the closures of the complementary intervals to points, then identifying the uncollapsed intervals in pairs, preserving their orientations, according to which pairs of points of $c\cap b$ are joined in the surgery process.  The coefficients $t_i$ define lengths on the edges of $T$, making it a metric tree.  The surgery process then continuously and canonically shrinks $T$ to a point by simultaneously shrinking all extremal outward-pointing edges at unit speed.  Each time an edge shrinks to a point, the resulting path $c_u(t)$ crosses a face of a simplex of $C_x(S)$.  From this viewpoint toward the surgery process it is clear that it yields continuous paths $c_u(t)$ that depend continuously on the coefficients $t_i$, including what happens when some $t_i$'s go to zero. 

In this way we obtain a deformation retraction of $C_x(S)$ onto the subcomplex spanned by cycles disjoint from $b$.  We claim this subcomplex is the star of $b$ in $C_x(S)$.  Certainly the subcomplex contains the star.  Conversely, let $[c_0,\cdots,c_k]$ be a simplex with each $c_i$ disjoint from $b$. Since $b$ is a single curve, it is contained in one of the cobordisms $C_i$ associated to $[c_0,\cdots,c_k]$.  We will show that $c_i$ and $b$ are cobordant by a cobordism $C\subset C_i$.  This will finish the proof of the claim, and hence of the theorem when $x$ is primitive, since it means that $[c_0,\cdots,c_i,b,c_{i+1},\cdots,c_k]$ is a simplex of $C_x(S)$.

The cycle $c_i$ represents the same homology class as $b$, so the difference $c_i - b$ is the boundary of a simplicial $2$-chain $\kappa$ in any triangulation of $S$ that contains $c_i$ and $b$ as subcomplexes.  We can think of $\kappa$ as a function assigning an integer to each complementary component of $c_i \cup b$, where this integer changes by $1$ when we cross $c_i$ or $b$.  The function $\kappa$ cannot take on more than two values, the two values on the complementary components adjacent to $b$, otherwise there would be a complementary region where $\kappa$ took on an extremal value different from the two values adjacent to $b$, and this region would have its boundary equal to a subcycle of $c_i$, an impossibility since $c_i$ is reduced. Since $\kappa$ takes on only two values, the union of the regions where it takes on a single value is a cobordism between $c_i$ and $b$.  If $k=0$ we are done, so we may assume $k>0$.  Let $C$ be the one of the two cobordisms between $c_i$ and $b$ that, near $c_i$, lies on the same side of $c_i$ as $C_i$.  Then we have $\bdy C = b - c_i$ and $\bdy C_i = c_{i+1} - c_i$, hence $\bdy(C - C_i) = b - c_{i+1}$.  We can write $C - C_i = A - B$ where $A$ is the closure of the set-theoretic difference $C\mathop{\backslash} C_i$ and $B$ is the closure of $C_i\mathop{\backslash} C$. Since $\bdy(A - B) = \bdy(C - C_i) = b - c_{i+1}$ and $b$ is part of the boundary of $B$ since it is contained in the interior of $C_i$, we must have $\bdy A \subset c_{i+1}$.  However, $c_{i+1}$ is reduced, so this implies $A$ is empty.  Thus $C\subset C_i$ and we are finished.   \qed

\smallskip
\pf Second proof of the Theorem, for $x$ nonzero.  It will suffice to show that each map $f\:S^n\to C_x(S)$ extends to a map $D^{n+1}\to C_x(S)$, for arbitrary $n\ge 0$.  We can assume $f$ is simplicial with respect to some triangulation of $S^n$.  First we want to choose a family $c(t)$ of weighted cycles representing the isotopy classes $f(t)$, having the form $c(t) =\sum_i t_i c_i(t)$ in each simplex of $S^n$ as in the first proof, but with each $c_i(t)$ now varying continuously as $t$ ranges over all of $S^n$, with the understanding that a cycle $c_i(t)$ is deleted from $c(t)$ when its coefficient $t_i$ becomes zero.  As a first approximation to constructing such a family $c(t)$ we can choose a hyperbolic structure on $S$ since we can assume its genus is greater than $1$, and then take the unique geodesic in each isotopy class of simple closed curves.  This does not quite work, however, since cycles can contain isotopic copies of the same curve, as can cycles of cycles.  To fix this, take a small $\epsilon$-neighborhood of each geodesic and choose appropriate parallel copies of the geodesic in this neighborhood.  All that matters here is the order of the parallel copies in the neighborhood, and this is specified by the isotopy classes of the cycles.  (The space of choices of parallel copies with a given ordering is contractible.)  

We can thicken each $c(t)$ to a bicollared neighborhood, varying continuously with $t$, and with the bicollar on a cycle $c_i(t)$ shrinking to thickness zero when the coefficient $t_i$ goes to zero.  Each bicollared $c(t)=\sum_i t_i C_i(t)$ determines a quotient map $g_t\:S\to S^1$  obtained by projecting the bicollar on each $c_i(t)$ to an interval of length $2\pi t_i$ and collapsing each complementary cobordism $C_i(t)$ to a point.  Since $\sum_i t_i = 1$, the quotient is a circle of radius $1$ with a preferred orientation, so it is identified with the standard $S^1$ up to a rotation.  The ambiguity of the rotation can be eliminated by fixing a basepoint in $S$ and letting $g_t$ send this to a chosen basepoint in $S^1$.  The map $g_t$ may not be smooth at the boundaries of the bicollars, but this can be achieved by reparametrizing the interval factors of the bicollars suitably.  We can assume the maps $g_t$ vary continuously in the $C^\infty$ topology on the space of smooth maps $S\to S^1$.  

By elementary obstruction theory, the space of basepoint-preserving maps $S\to S^1$ has path-components in bijective correspondence with $H^1(S)\approx H_1(S)$ and all its homotopy groups $\pi_n$ for $n>0$ trivial.  This implies that we can extend the family $g_t$ to a family $g_t\:S\to S^1$ for $t\in D^{n+1}$.  These $g_t$'s can be taken to be $C^\infty$ maps, varying continuously in the $C^\infty$ topology.

We can construct a family $c(t)$ over $D^{n+1}$ by taking linear combinations of preimages of regular values of $g_t$ in the following way.  For each $t\in D^{n+1}$ the regular values of $g_t$ are dense by Sard's theorem, so we can choose a regular value $v_t$.  This remains a regular value in some open neighborhood $U_t$ of $t$ in $D^{n+1}$.  Doing this for all $t$, we obtain an open cover of $D^{n+1}$ by the sets $U_t$.  Since $D^{n+1}$ is compact, a finite number of $U_t$'s suffice to cover $D^{n+1}$.  We relabel these as $U_j$.  Thus for each $j$ we have a regular value $v_j$ and a family of cycles $c_j(t)=g_t^{-1}(v_j)$ in the homology class $x$, for $t\in U_j$.  The cycles $c_j(t)$ for all $U_j$ containing a fixed $t$ then give a cycle of cycles since the preimages of the intervals between adjacent $v_j$'s in $S^1$ give the necessary cobordisms between these cycles $c_j(t)$.  Choosing a partition of unity subordinate to this cover, consisting of functions $\phi_j$ with support in $U_j$ and $\sum_j\phi_j=1$, we can then form the sums $c(t) =\sum_j\phi_j(t) c_j(t)$.  If all the cycles $c_j(t)$ were reduced we would then have a family of points $c(t)$ in $C_x(S)$ varying continuously with $t\in D^{n+1}$.

Let us describe a reduction procedure for making each $c(t)$ into a cycle of reduced cycles.  Consider the dual graph $\Gamma_t$ of $c(t)$.  This has edges corresponding to the components of the $c_j(t)$'s and vertices corresponding to the components of the complementary cobordisms.  The edges are oriented according to the normal orientations of the $c_j(t)$'s and the coefficients $\phi_j(t)$ assign lengths to the edges, so $\Gamma_t$ is an oriented metric graph. Vertices of $\Gamma_t$ that are sources or sinks correspond to complementary components of $c(t)$ that are null cobordisms.  We first eliminate the sinks by shrinking all edges leading to them at unit speed until these edges collapse to points, then continuing on with any remaining sinks, eventually reaching a new $\Gamma_t$ with no sinks.  When the edges surrounding a sink reach zero length we delete the corresponding components of $c(t)$ and amalgamate the complementary component corresponding to the sink vertex with the components corresponding to the adjacent vertices.  After all sinks are eliminated we have a new cycle of cycles, and we do the analogous procedure on it to eliminate all source vertices.  This does not create new sinks.  In the end we have a new cycle of reduced cycles which we still call $c(t)$.  It is evident that the reduction process depends continuously on $t$.  During the reduction process the functions $g_t$ can be deformed each time a sink or source is eliminated, so that the new family $c(t)$ consists of preimages under the new $g_t$ of the same regular values $v_j$.  Then for $t\in S^n$ the new $c(t)$ and the original family $c(t)=f(t)$ both consist of preimages of regular values of $g(t)$ for all $t$, so these two families are homotopic in $C_x(S)$ by the obvious linear homotopy. Thus the original map $f\:S^n\to C_x(S)$ is nullhomotopic \qed

\section References 

\parindent 45pt

\item{[BBM]} M. Bestvina, K.-U. Bux, D. Margalit, The dimension of the Torelli group. arXiv:0709.0287.

\item{[H]} A. Hatcher, On triangulations of surfaces. Topology Appl. 40 (1991) 189-194. Revision: http://www.math.cornell.edu/$\tilde{\phantom{x}}$hatcher/Papers/TriangSurf.pdf

\end